\documentclass[a4paper,12pt]{article}
\usepackage{amssymb}

\newcommand{\R}{\mathbb{R}}

\newcommand{\cuad}{{\sqcap\kern-.68em\sqcup}}

\newcommand{\norm}[1]{\|#1\|}

\newtheorem{theorem}{Theorem}[section]
\newtheorem{proposition}{Proposition}[section]

\newtheorem{lemma}{Lemma}[section]
\newtheorem{corollary}{Corollary}[section]
\newtheorem{remark}{Remark}[section]
\newcommand{\bremark}{\begin{remark} \em}
\newcommand{\eremark}{\end{remark} }

\begin{document}

\begin{center}{\bf  \large   Fundamental solutions for Schr\"{o}dinger operators with     \\[2mm]

general inverse square potentials  }\medskip
  \bigskip\medskip

{\small

 Huyuan Chen\footnote{chenhuyuan@yeah.net} \ \,
  Suad Alhomedan\footnote{shemedan@ksu.edu.sa}\ \,   Hichem Hajaiej\footnote{hichem.hajaiej@gmail.com, hichem.hajaiej@nyu.edu} \ \,    Peter Markowich\footnote{ peter.markowich@kaust.edu.sa}

\medskip \medskip

$^{1}$ Department of Mathematics, Jiangxi Normal University,  Nanchang,\\
 Jiangxi 330022, PR China

 \smallskip

$^{2}$    Suad Alhemedan, Department of Mathematics, College of Science,\\ King Saud University

 \smallskip

$^{3}$  New York University  Shanghai 1555  Century Avenue Pudong New Area,\\
Shanghai, China

 \smallskip

$^{4}$    Division of Math and Computer Sci. and Eng., King Abdullah University\\ of Science
and Technology, Thuwal 23955-6900, Saudi Arabia
}


\begin{abstract}
In this paper,  we clarify the fundamental solutions for  Schr\"{o}dinger operators   given as
$\mathcal{L}_{\mu V}=-\Delta- \mu  V$, where the potential $V$ is   a general inverse square potential in $\R^N\setminus\{0\}$ with $N\ge 3$.
In particular, letting $\mu\in(0,\, (N-2)^2/4]$,
$$\lim_{x\to0}V(x)|x|^{2}=1 \quad {\rm and}\quad \lim_{|x|\to+\infty} V(x)|x|^2=t,$$
where $t\ge0$, we discuss the existence and nonexistence of positive fundamental solutions for Hardy operator $\mathcal{L}_{\mu V}$, which depend  on the parameter $t$.

\end{abstract}

\end{center}

  \vspace{1mm}


\vspace{2mm}

\setcounter{equation}{0}
\section{Introduction}

In this paper, we are concerned with the fundamental solutions for   \textbf{Schr\"{o}dinger operators}  given as
\begin{equation}\label{eq 1.0}
 \mathcal{L}_{\mu V}=-\Delta-  \mu V,
\end{equation}
where $0<\mu\le \mu_0:=\frac{(N-2)^2}{4}$ with $N\ge 3$ and $V\in C^\gamma_{loc}(\R^N\setminus\{0\})$  with $\gamma\in(0,1)$ is a nonnegative \textbf{general inverse square  potential}.

When  $\mu\le \mu_0$ and $V(x)=|x|^{-2}$ for $x\not=0$, we reduce the notation $\mathcal{L}_{\mu V}$  into $\mathcal{L}_{\mu}$, given as,
$$\mathcal{L}_\mu:= -\Delta   -\frac{\mu}{|x|^{2 }}.$$
 The corresponding   Hardy inequalities, see \cite{ACR,BM,DD,GP},  enable the variational techniques available to solve
the semilinear elliptic and parabolic  differential equations,  for instance\cite{BOP,CC,D,FF,KV,IYM,MT}. In bounded smooth domain,  Comparison Principle for $\mathcal{L}_{\mu}$
holds by the Hardy inequalities; but it fails in unbounded domain. From the observation of fundamental solutions for $\mathcal{L}_\mu$,  the authors in \cite{BDT}  studied the existence and nonexistence of  isolated singular solutions for semilinear Hardy problem
$\mathcal{L}_{\mu } u=u^p\ \, {\rm in}\ \, B_1(0)\setminus\{0\}.$

 It is known that when $\mu\le \mu_0$,  the problem
 \begin{equation}\label{eq 1.2}
\mathcal{L}_\mu u= 0\quad {\rm in}\quad  \R^N\setminus \{0\}
\end{equation}
has two  radially symmetric solutions with the explicit formulas that
\begin{equation}\label{f1.1}
 \Phi_\mu(x)=\left\{\arraycolsep=1pt
\begin{array}{lll}
 |x|^{\tau_-(\mu)}\quad
   &{\rm if}\quad \mu<\mu_0\\[1mm]
 \phantom{   }
|x|^{\tau_-(\mu)}(-\ln|x|) \quad  &{\rm   if}\quad \mu=\mu_0
 \end{array}
 \right.\qquad {\rm and}\qquad \Gamma_\mu(x)=|x|^{ \tau_+(\mu)},
\end{equation}
where
$$
 \tau_-(\mu)=-\frac{N-2}2-\sqrt{\mu_0-\mu}\quad{\rm and}\quad   \tau_+(\mu)=-\frac{ N-2}2+\sqrt{\mu_0-\mu}.
$$
Here the parameters $\tau_-(\mu)$ and $ \tau_+(\mu)$ are the zero points of
$\tau(\tau+N-2)+\mu=0.$

Recently,  the authors in \cite[Theorem  1.1]{CQZ} showed that the  solution $\Phi_\mu$ verifies the  $d\mu$-distributional identity
$$
 \int_{\R^N}\Phi_\mu   \mathcal{L}^*_\mu(\xi) \, d\mu  =c_\mu\xi(0),\quad \forall\, \xi\in C^{1.1}_c(\R^N),
$$
where  $d\mu(x) =\Gamma_\mu(x) dx$,
\begin{equation}\label{L}
 \mathcal{L}^*_\mu=-\Delta -2\frac{\tau_+(\mu) }{|x|^2}\,x\cdot\nabla
\end{equation}
and
 \begin{equation}\label{cmu}
 c_\mu=\left\{\arraycolsep=1pt
\begin{array}{lll}
 2\sqrt{\mu_0-\mu}\,|\mathcal{S}^{N-1}|\quad
   &{\rm if}\quad \mu<\mu_0,\\[2mm]
 \phantom{   }
|\mathcal{S}^{N-1}| \quad  &{\rm  if}\quad \mu=\mu_0.
 \end{array}
 \right.
 \end{equation}
 Here  $\mathcal{S}^{N-1}$ is the sphere of  the unit ball of $\R^N$ and $|\mathcal{S}^{N-1}|$ is the volume of the unit sphere.
 Normally, $\Phi_\mu$ is viewed as a fundamental solution of $\mathcal{L}_\mu$. We note that the fundamental solution $\Phi_\mu$ keeps positive
when $\mu<\mu_0$ and changes signs for $\mu=\mu_0$.

Our purpose of this article is to consider the fundamental solutions for   $\mathcal{L}_{\mu V}$. To this end,
 we have to investigate the singular  solutions of
\begin{equation}\label{eq 1.1}
    \mathcal{L}_{\mu V}\, u=0 \quad    {\rm in}\quad \R^N\setminus\{0\}
\end{equation}
under the hypotheses that
$0<\mu\le \mu_0$ and $V\in C^\gamma_{loc}(\R^N\setminus\{0\})$   verifies that
\begin{equation}\label{V0}
V(x)\ge |x|^{-2}\quad{\rm and}\quad \lim_{x\to0}V(x)|x|^{2}=1.
\end{equation}
 Note that  the operator $\mathcal{L}_{\mu V}$ could be viewed as a perturbation of $\mathcal{L}_\mu$ near the origin, that is,
$$\mathcal{L}_{\mu V}=\mathcal{L}_\mu-\mu(V-|x|^{-2}). $$

The following proposition is to clarify the classical solution to (\ref{eq 1.1}) in the $d\mu$-distributional sense.

\begin{proposition}\label{classification}
Assume that  $N\ge3$, $ 0< \mu\le \mu_0$ and $V\in C^\gamma_{loc}(\R^N\setminus\{0\})$ is a positive Hardy potential satisfying $(\ref{V0})$ and $u$ is a classical solution of (\ref{eq 1.1}) such that
$$u\ge0\quad{\rm in}\quad B_r(0)\setminus\{0\}$$
for some $r>0$.
Then
\begin{equation}\label{v0}
 \int_{B_1(0)}\left(V(x)-|x|^{-2}\right)|x|^{2-N} dx<+\infty
\end{equation}
and
 there exists $k\ge0$ such that $u_k$ is a distributional solution of
 \begin{equation}\label{eq 1.1w}
  \Gamma_\mu\,  \mathcal{L}_{\mu V}\, u_k=k\delta_0 \quad {\rm in  } \ \ \R^N,
\end{equation}
i.e.
\begin{equation}\label{d 1.2}
 \int_{\R^N} u_k \mathcal{L}^*_{\mu V}(\xi) \, d\mu  =   c_\mu k\xi(0),\quad \forall\, \xi\in C^{1.1}_c(\R^N),
\end{equation}
 where $d\mu(x) =\Gamma_\mu(x)\, dx$,
 $$\mathcal{L}^*_{\mu V}=\mathcal{L}^*_{\mu}-\mu (V-|x|^{-2})$$
 and $c_\mu$ is given in   (\ref{cmu}).

 Furthermore, if $\mu\in(0,\mu_0)$, $k>0$ and the $d\mu$-distributional solution $u_k\ge0$ in $\R^N\setminus\{0\}$, then
  $$u_k\ge k\Phi_\mu\quad {\rm in}\quad \R^N\setminus\{0\}.$$
\end{proposition}

  Our concern is  the fundamental solutions of  \textbf{ Schr\"{o}dinger operator} $\mathcal{L}_{\mu V}$.  Motivated by Proposition \ref{classification}, we introduce the notation of   fundamental solutions of $\mathcal{L}_{\mu V}$  with $\mu\le \mu_0$ as following: {\it A function $u$ is said to be a fundamental solution of $\mathcal{L}_{\mu V}$, if $u$ is a classical solution of (\ref{eq 1.1}) and verifies the  identity (\ref{d 1.2})
 with  $k\not=0$. Furthermore, $u$ is said to be a positive fundamental solution of  $\mathcal{L}_{\mu V}$,  if $u>0$ is a fundamental solution with $k>0$ in (\ref{d 1.2}).  }
In \cite{CQZ}, it shows that $\Phi_\mu$ with $\mu\in(-\infty,\,\mu_0)$ is a positive fundamental solution of $\mathcal{L}_\mu$ verifying (\ref{d 1.2}) with  $k=1$, and $\Phi_{\mu_0}$ is a signs-changing fundamental solution of  $\mathcal{L}_{\mu_0}$ with $k=1$.

We first consider   fundamental solutions of (\ref{eq 1.1}), when  $V(x)|x|^{2}$ is not too large  at infinity. To be precise, we propose the following assumption:
\begin{itemize}
\item[$(V)$]
the Hardy potential $V\in C^\gamma_{loc}(\R^N\setminus\{0\})$ with $\gamma\in(0,1)$ satisfies that
\begin{enumerate}\item[$(i)$] there exists $c_0>0$ such that
$$   0\le  V(x)-|x|^{-2}  \le  c_0,\quad\ \forall\, x\in \R^N\setminus\{0\};$$
\end{enumerate}
\begin{enumerate}\item[$(ii)$] for some $\alpha\in[1,\frac{\mu_0}{\mu})$,
$$\limsup_{|x|\to+\infty} V(x)|x|^2=\alpha.$$
\end{enumerate}
\end{itemize}

  The existence of fundamental solution of $\mathcal{L}_{\mu V}$ states as following.

\begin{theorem}\label{teo 2}
Assume that   $N\ge3$, $ 0< \mu<\mu_0$   and $V$ satisfies the assumption $(V)$.
Then  $\mathcal{L}_{\mu V}$
 has positive  fundamental solutions and for  $k>0$, there is a minimal positive solution $u_k$ verifying (\ref{d 1.2}) with such $k$. Moreover,
 \begin{equation}\label{1.3}
\lim_{x\to0}u_k(x)\Phi_\mu^{-1}(x)=k
 \end{equation}
 and for any $\mu'\in (\alpha\mu,\, \mu_0)$, there exists  $c_1>0$ such that
  \begin{equation}\label{1.3.1}
 k\Phi_\mu\le u_k\le c_1\Phi_{\mu'}\quad {\rm in}\quad \R^N\setminus B_1(0).
 \end{equation}
\end{theorem}

Here   $u_k$ is said to be a minimal positive fundamental solution of $\mathcal{L}_{\mu V}$ verifying (\ref{d 1.2}) with  $k>0$ if any positive fundamental solution $u$ of $\mathcal{L}_{\mu V}$ verifying (\ref{d 1.2}) with such $k$, verifies that $u\ge u_k$.

We note that for $0<\mu<\mu_0$, $k>0$ fixed, we can not get the uniqueness of the fundamental solutions of $\mathcal{L}_{\mu V}$ verifying (\ref{d 1.2}) with such $k$, because of the lack of  weak comparison principle in unbounded domain.  In particular,   $\{\Phi_\mu+l\Gamma_\mu\}$ with $l\ge0$ is the set of positive fundamental solutions for  $\mathcal{L}_{\mu}$ verifying (\ref{d 1.2}) with  $k=1$, even subjecting to $\lim_{|x|\to+\infty}u(x)=0$, and $\Phi_\mu$ is the minimal positive fundamental solution verifying (\ref{d 1.2}) with  $k=1$.

Involving general inverse square potential $V$ satisfying $(V)$,  the minimal positive fundamental solution verifying (\ref{d 1.2}) with $k>0$ is derived by iterating the sequence
$$w_0=k\mathcal{G}_\mu\quad{\rm and}\quad w_n=\mu \mathbb{G}_{B_n(0)}[Vw_{n-1}],$$
where $\mathcal{G}_\mu=\Phi_\mu-\Gamma_\mu$,  the fundamental solution of $\mathcal{L}_{\mu V}$ in $B_1(0)$ with Dirichlet boundary condition and  $\mathbb{G}_{B_n(0)}$ is the Green operator defined by the Green kernel $G_{B_n(0)}$.
 A super bound is constructed to control this sequence by using the hypothesis that $\alpha<\frac{\mu_0}{\mu}$ in $(V)$. This super bound also provides estimates for the singularity at the origin and the decay  at infinity of the minimal fundamental solutions. The isolated singularity  (\ref{1.3}) at origin is also motivated  by the classification of isolated singularities in \cite[Proposition 4.2]{CQZ}. \smallskip

On the contrary, if $V(x)|x|^2$ is large enough at infinity,\textbf{ how  is it going on}   the fundamental solutions of the operator $\mathcal{L}_{\mu V}$?
From \cite{CQZ}, it shows that if $\mu V(x)|x|^2>\mu_0$ near the origin, there exist no fundamental solutions for $\mathcal{L}_{\mu V}$.
Our second purpose is to \textbf{study} the nonexistence of fundamental solutions of $\mathcal{L}_{\mu V}$ when  $V(x)|x|^{2}$ is large  at infinity.  Precisely, we assume that
\begin{itemize}
\item[$(\widetilde{V})$]
the Hardy potential  $V\in C^\gamma_{loc}(\R^N\setminus\{0\})$ with $\gamma\in(0,1)$ satisfies that
\begin{enumerate}\item[$(I)$] there exists $c_0>0$   such that
$$  0 \le V(x)-|x|^{-2}\le c_0,\quad\ \forall\, x\in \R^N\setminus\{0\};$$
\end{enumerate}
\begin{enumerate}\item[$(II)$] for some $\beta>1$,
$$\liminf_{|x|\to+\infty} V(x)|x|^2=\beta.$$
\end{enumerate}
\end{itemize}

Our maim result on the nonexistence states as following.

\begin{theorem}\label{teo 1}
Assume that $N\ge3$, $ 0< \mu<\mu_0$ and $V$   satisfies $(\widetilde{V})$.
Then  there exists $\beta^*>1$ such that for $\beta>\beta^*$, the Hardy operator $\mathcal{L}_{\mu V}$
 has no positive fundamental solutions.

\end{theorem}

We note that the nonexistence result is derived by considering the decay at infinity. By contradiction, if there is a  positive fundamental solution $u$ for $\mathcal{L}_{\mu V}$, our strategy is to iterate  an origin decay $k\Phi_\mu $,   then  to improve the coefficient of the decay.  In fact, for $0<\mu<\mu_0$ and $\liminf_{|x|\to+\infty} V(x)|x|^2>1$, the positive fundamental solution of (\ref{eq 1.1})
has the asymptotic  behavior that
\begin{equation}\label{1.1}
 \lim_{|x|\to+\infty} u(x)|x|^{-\tau_-(\mu)}=+\infty.
\end{equation}
Finally,  a contradiction is deduced  by the choice of $\beta^*$.

To make clear of the connection between Theorem \ref{teo 2} and Theorem \ref{teo 1}, we have build the following results on the fundamental solution for
  a typical  class of Hardy operators.

\begin{corollary}\label{cr 1.1}
Assume that   $N\ge3$, $ 0< \mu<\mu_0$, $\rho\ge1$   and
\begin{equation}\label{Vr}
 V_\rho(x)=\frac1{|x|^2}\,\frac{1+\rho|x|^2}{1+|x|^2},\qquad \forall\, x\in\R^N\setminus\{0\}.
\end{equation}
Then   there exists $\rho^*\ge \frac{\mu_0}{\mu}$ such that for $\rho<\rho^*$, $\mathcal{L}_{\mu V_\rho}$
 has positive fundamental solutions and  for  $\rho>\rho^*$, $\mathcal{L}_{\mu V_\rho}$ admits no positive fundamental solutions.

\end{corollary}

 \begin{remark}
 We note that $\{V_\rho\}_\rho$ is an increasing potential satisfying
  $$\lim_{x\to0} V_\rho(x)|x|^{2}=1\quad{\rm and}\quad \lim_{|x|\to+\infty}V_\rho(x)|x|^2=\rho.$$
 It is an interesting but open question whether $\rho^*=\frac{\mu_0}{\mu}$.
 \end{remark}

Furthermore, we discuss the existence of positive fundamental solutions for  \textbf{Schr\"{o}dinger operators} $\mathcal{L}_{\mu V}$ in the case that $V(x)\le |x|^{-2}$ and $\mu\in(0,\mu_0]$ in Section \S5. In particular, when $\mu=\mu_0$, Hardy operator $\mathcal{L}_{\mu_0}$ has a signs-changing  fundamental solution $\Phi_{\mu_0}$, so our concentration is on the existence of positive fundamental solutions when the potential
$\limsup_{|x|\to+\infty}V(x)|x|^{2}<1$.\smallskip

The rest of the paper is organized as follows. In Section 2,  we  clarify the classical solution of $\mathcal{L}_\mu u=f$ in the $d\mu$-distributional sense and prove Proposition \ref{classification}.   Section  3 is devoted
to prove the existence of fundamental solutions of $\mathcal{L}_{\mu V}$ under the assumption of $(V)$ and  the nonexistence of fundamental solutions
 of $\mathcal{L}_{\mu V}$  when $V$ verifies $(\widetilde{V})$. In Section \S4, we discuss the fundamental solutions of\textbf{ Schr\"{o}dinger operator} when   $V(x)\le |x|^{-2}$
 and the case $\mu=\mu_0$. Finally, we put the classification of isolated singularities at the origin of $\mathcal{L}_\mu u=f$ in $\R^N\setminus\{0\}$ in the Appendix.

\setcounter{equation}{0}
\section{Preliminary  }

Our aim in this section is  to prove Proposition \ref{classification}.  We first collect comparison principles in bounded domain. The first version states as following.
\begin{lemma}\label{lm ocp}\cite[Lemma 2.1]{CQZ}
Let $\mu\le\mu_0$, $O$ be a bounded open set in $\R^N$, $L: O\times [0,+\infty)\to[0,+\infty)$ be a continuous function satisfying that for any $x\in  O$,
$$L(x,s_1)\ge L(x,s_2)\quad {\rm if}\quad s_1\ge s_2,$$ then $\mathcal{L}_\mu+L$ with $\mu\ge \mu_0$ verifies the Comparison Principle,
that is, if
$$u,\,v\in C^{1,1}(O)\cap C(\bar O)$$ verify that
$$\mathcal{L}_\mu u+ L(x,u)\ge \mathcal{L}_\mu v+ L(x,v) \quad {\rm in}\quad  O
\qquad{\rm and}\qquad   u\ge  v\quad {\rm on}\quad \partial O,$$
then
$$u\ge v\quad{\rm in}\quad  O.$$

\end{lemma}

\begin{lemma}\label{lm cp}
Assume that $\mu\le\mu_0$,  $\Omega$ is a bounded $C^2$ domain containing the origin,     $f_1$, $f_2$ are two functions in $C^\gamma (\Omega\setminus\{0\})$ with $\gamma\in(0,1)$, $g_1$, $g_2$ are two continuous  functions on $\partial \Omega$, and
$$ f_1\ge f_2\quad {\rm in}\quad \Omega\setminus\{0\} \quad{\rm and}\quad  g_1\ge g_2\quad {\rm on}\quad \partial \Omega.$$
Let $u_i$ with $i=1,2$  be the classical solutions of
$$
 \arraycolsep=1pt\left\{
\begin{array}{lll}
 \displaystyle \mathcal{L}_\mu u = f_i\qquad
   {\rm in}\quad  {\Omega}\setminus \{0\},\\[1.5mm]
 \phantom{ L_\mu     }
 \displaystyle  u= g_i\qquad  {\rm   on}\quad \partial{\Omega}.
 \end{array}\right.
$$

If $$\liminf_{x\to0}u_1(x)\Phi_\mu^{-1}(x) \ge \limsup_{x\to0}u_2(x)\Phi_\mu^{-1}(x),$$
then
$$u_1\ge u_2\quad {\rm in}\quad \Omega\setminus\{0\}.$$

\end{lemma}
{\bf Proof.}  Let $w=u_2-u_1$ be a solution of
 $$\arraycolsep=1pt\left\{
\begin{array}{lll}
 \displaystyle \mathcal{L}_\mu u \le 0\qquad
   {\rm in}\quad  {\Omega}\setminus \{0\},\\[1.5mm]
 \phantom{ L_\mu     }
 \displaystyle  u\le 0 \qquad  {\rm   on}\quad \partial{\Omega},\\[1.5mm]
 \phantom{   }
  \displaystyle \limsup_{x\to0}u(x)\Phi_\mu^{-1}(x)\le0,
 \end{array}\right.$$
  then for any $\epsilon>0$, there exists $r_\epsilon>0$ converging to zero as $\epsilon\to0$ such that
 $$w\le \epsilon \Phi_\mu\quad{\rm on}\quad \partial B_{r_\epsilon}(0).$$
We see that
$$w\le 0<\epsilon \Phi_\mu \quad{\rm on}\quad \partial\Omega,$$
then by Lemma \ref{lm ocp}, we have that
$$w\le \epsilon \Phi_\mu\quad{\rm in}\quad \Omega\setminus\{0\}. $$
By the arbitrary of $\epsilon>0$, we have that $w\le 0$   in $\Omega\setminus\{0\}$.
\hfill$\Box$\medskip

It is remarkable that when $\Omega=\R^N\setminus\{0\}$, the comparison principle for $\mathcal{L}_\mu$ fails for $\mu>0$, even subject to the condition that
$\lim_{|x|\to+\infty}u(x)=0.$  A counterexample is that  $\Phi_\mu$ and $\Gamma_\mu$ are the   solutions of $\mathcal{L}_\mu u=0$ in $\R^N\setminus\{0\}$.

Motivated by the Kato's inequality, see Proposition 6.1 in \cite{P}, we have the following  comparison principle in the weak sense.

\begin{lemma}\label{lm wcp}
Assume that $\mu\le\mu_0$,   $\Omega$ is a bounded $C^2$ domain containing the origin,  $k_1\ge k_2$,  $f_1$, $f_2$ are two functions in $C^\gamma (\Omega\setminus\{0\})$ with $\gamma\in(0,1)$, $g_1$, $g_2$ are two continuous  functions on $\partial \Omega$, and
$$ f_1\ge f_2\quad {\rm in}\quad \Omega\setminus\{0\} \quad{\rm and}\quad  g_1\ge g_2\quad {\rm on}\quad \partial \Omega.$$
Let $u_i$ with $i=1,2$ be the $d\mu$-distributional solutions  of
\begin{equation}\label{eq1 2.1}
 \arraycolsep=1pt\left\{
\begin{array}{lll}
 \displaystyle \mathcal{L}_\mu u = f_i+k_i\delta_0\qquad
   {\rm in}\quad  {\Omega},\\[1.5mm]
 \phantom{ L_\mu     }
 \displaystyle  u= g_i\qquad  {\rm   on}\quad \partial{\Omega}.
 \end{array}\right.
\end{equation}
i.e.
 $$\int_\Omega u_i\mathcal{L}^*_\mu \xi\, d\mu=\int_{\Omega} f_i \xi \,d\mu +\int_{\partial\Omega}g_i \frac{\partial  \xi }{\partial\nu}\Gamma_\mu d\omega +k_i\xi(0),\quad\forall\, \xi\in C^{1,1}_0(\Omega). $$

Then
$$u_1\ge u_2\quad {\rm in}\quad \Omega\setminus\{0\}.$$

\end{lemma}
{\bf Proof.}    Let $w=u_2-u_1$, then we have that
  $$\int_\Omega w\mathcal{L}^*_\mu \xi\, d\mu=\int_{\Omega} (f_2-f_1)\xi \,d\mu +\int_{\partial\Omega}(g_2-g_1)\frac{\partial \xi}{\partial\nu}\Gamma_\mu\, d\omega +(k_2-k_1)\xi(0),\quad\forall\, \xi\in C^{1,1}_0(\Omega), $$
where $\nu$ is the unit normal vector pointing outside of $\Omega$. So we now put  $\xi\in C^{1,1}_0(\Omega)$, $\xi\ge0$,  then  we have that $\frac{\partial  \xi }{\partial\nu}\ge0$
on $\partial\Omega$ and
$$\int_\Omega w\mathcal{L}^*_\mu \xi\, d\mu\le 0.$$

 Taking  $O=\{x\in\Omega:\, w(x)>0\}$, denote by $\eta_{\omega,n}$ the
solution of
\begin{equation}\label{L00}
\arraycolsep=1pt\left\{
\begin{array}{lll}
 \displaystyle  \mathcal{L}_\mu^*  u= \zeta_n\qquad
   &{\rm in}\quad   \Omega,\\[1mm]
 \phantom{  \mathcal{L}_\mu^*  }
 \displaystyle  u= 0\quad & {\rm   on}\quad \partial{\Omega},
 \end{array}\right.
 \end{equation}
where $\zeta_n:\bar\Omega\mapsto[0,1]$ is a nonnegative, $C^1(\bar\Omega)$
function such that
$$\zeta_n\to\chi_O\quad {{\rm in}}\ L^\infty( \Omega)\quad {{\rm as}}\  \ n\to\infty.$$
Then
$$\displaystyle\int_{\Omega}w\Gamma_\mu\, \zeta_n \ dx\le 0.
$$
Then passing to the limit as $n\to\infty$, we have
$$\displaystyle\int_{O}w\Gamma_\mu \ dx\le 0.$$ This implies $w\le 0$ a.e. in $\Omega$.
The proof ends.\hfill$\Box$\medskip

To estimate the isolated singularities at the origin of the fundamental solutions, we   clarify the isolated singular solutions of
\begin{equation}\label{eq 2.1.1}
  \arraycolsep=1pt\left\{
\begin{array}{lll}
 \displaystyle   \ \ \mathcal{L}_\mu\, u=f \quad
   {\rm in}\quad \R^N\setminus\{0\},\\[2mm]
 \phantom{    }
 \displaystyle \lim_{|x|\to+\infty}u(x)=0.
 \end{array}\right.
\end{equation}

\begin{lemma}\label{pr 2.1}
Assume that $N\ge 2$, $ \mu\le \mu_0$, $f$ is a function in $C^\gamma_{loc}(\R^N\setminus \{0\})\cap L^1_{loc}(\R^N,\, d\mu(x))$
and
$u$ is a  solution of (\ref{eq 2.1.1}) satisfying that $u\ge0$ in $B_r(0)\setminus\{0\}$ for some $r>0$.
  Then there exists $k\ge 0$ such that  $u$ verifies the $d\mu$-distributional identity
   \begin{equation}\label{1.2f}
 \int_{\R^N}u  \mathcal{L}_\mu^*\xi\, d\mu  = \int_{\R^N} f  \xi\, d\mu +c_\mu k\xi(0),\quad\forall\, \xi\in   C^{1.1}_0(\R^N).
\end{equation}
\end{lemma}
{\bf Proof.}  In a bounded domain, this proposition could be seen in  Proposition 4.2 in \cite{CQZ}. For the convenience, we give the details in the appendix.
\hfill$\Box$\medskip

\begin{lemma}\label{lm 2.3}
Assume that $N\ge 3$, $-\infty<\mu<\mu_0$, $f\ge0$ and $u_k$ is a nonnegative solution of (\ref{eq 2.1.1}) and verifies (\ref{1.2f}) with  $k\ge 0$. Then
\begin{equation}\label{2.6}
 u_k\ge k\Phi_\mu\quad {\rm in }\quad \R^N\setminus\{0\}.
\end{equation}
\end{lemma}
{\bf Proof.} Let $$\mathcal{G}_\mu=\Phi_\mu-\Gamma_\mu\quad {\rm in}\quad B_1(0),$$
which is
  the $d\mu$-distributional solution of
  \begin{equation}\label{eq 2.2}
 \arraycolsep=1pt\left\{
\begin{array}{lll}
 \displaystyle  \mathcal{L}_\mu u =c_\mu \delta_0\qquad
   &{\rm in}\quad   B_1(0),\\[1.5mm]
 \phantom{ \mathcal{L}_\mu    }
 \displaystyle  u= 0\qquad  &{\rm   on}\quad \partial{ B_1(0)}.
 \end{array}\right.
 \end{equation}

 Denote by
 $$\mathcal{G}_{\mu,r}(x)=r^{\tau_-(\mu)}\mathcal{G}_\mu(r^{-1}x)=|x|^{\tau_-(\mu)}-r^{\tau_-(\mu)-\tau_+(\mu)}|x|^{\tau_+(\mu)}, $$
then $\mathcal{G}_{\mu,r}$ is the solution of (\ref{eq 2.2}) replacing $B_1(0)$ by $B_r(0)$, and by Lemma \ref{lm wcp} for any $r>1$, we have that
\begin{equation}\label{2.1.1}
 u(x)\ge k\mathcal{G}_{\mu,r}(x),\quad\ \forall\, x\in B_r(0).
\end{equation}
By the fact that
$$\lim_{r\to+\infty}\mathcal{G}_{\mu,r}(x)= |x|^{\tau_-(\mu)}, $$
passing to the limit of (\ref{2.1.1}) as $r\to+\infty$, it
implies (\ref{2.6}).\hfill$\Box$\medskip

For $N\ge3$, it is known that $\Phi_0(x)=c_N|x|^{2-N}$
is the fundamental solution of
$$-\Delta \Phi_0=\delta_0\quad{\rm in}\quad \mathcal{D}'(\R^N), $$
where $c_N>0$ is the normalized constant. Let $G_\Omega$ be the Green kernel of $-\Delta$ in $\Omega\times\Omega$, subjecting to the Dirichlet boundary condition.
We define by $\mathbb{G}_\Omega$   the Green operator defined by Green kernel $G_\Omega$
$$\mathbb{G}_\Omega[f](x)=\int_{\Omega}G_\Omega(x,y) f(y)dy, \quad \forall f\in C^\infty_c(\Omega).$$
In particular, we have that
$$G_{\R^N}(x,y)= c_N|x-y|^{2-N},\qquad\forall\, x,y\in\R^N,\, x\not=y.$$

\begin{lemma}\label{re 2.1}
Assume that $N\ge3$,  $u$ is a nonnegative classical solution of
\begin{equation}\label{eq}
 -\Delta u=Qu\quad {\rm in}\quad \R^N\setminus\{0\},
\end{equation}
where $Q>0$ is locally H\"{o}lder continuous in $\R^N\setminus\{0\}$.
If
\begin{equation}\label{e 2.1}
 \lim_{r\to0^+}\int_{B_1(0)\setminus B_r(0)} Q(x)|x|^{2-N} dx=+\infty,
\end{equation}
then
 $$u=\mathbb{G}_{\R^N}[Qu].$$
\end{lemma}
{\bf Proof.} We first claim that $\mathbb{G}_{\R^N}[Qu]$ is well-defined in $\R^N\setminus\{0\}$, that is,
\begin{equation}\label{2.04}
 \int_{\R^N}\frac{Q(y)u(y)}{|x-y|^{N-2}}\,dy<+\infty,\quad\forall\, x\in \R^N\setminus\{0\}.
\end{equation}
At this moment, we assume that the claim  (\ref{2.04})  holds. Fix a point $x_0\in \partial B_1(0)$, then (\ref{2.04}) implies that
\begin{equation}\label{2.05}
 \int_{B_{2}(0)} Q(y)u(y)\, dy<+\infty,
\end{equation}
so letting $f:=Qu\in L^1_{loc}(\R^N)$, applying Theorem 1.1 in \cite{BL}, there exists $k\ge0$ such that
$$-\Delta u=Qu+k\delta_0\quad{\rm in}\quad \mathcal{D}'(\R^N) $$
and if $k>0$, by Lemma \ref{lm wcp}, we have that
$$u\ge k \mathbb{G}_{B_2(0)}[\delta_0].$$
Observing that
$$\mathbb{G}_{B_2(0)}[\delta_0](x)\ge c_2|x|^{2-N},\quad\forall\, x\in B_1(0)\setminus\{0\},$$
then
\begin{eqnarray*}
\int_{B_1(0)\setminus B_{r}(0)}Q(y)u(y)\, dy &\ge &c_2 k \int_{B_1(0)\setminus B_{r}(0)}Q(y)|y|^{2-N}\, dy \\
   &\to& +\infty\quad{\rm as}\quad r\to0^+,
\end{eqnarray*}
 which contradicts (\ref{2.05}). So we have that $k=0$, that is,
 $$-\Delta u=Qu \quad{\rm in}\quad \mathcal{D}'(\R^N). $$
So we have that
$$u=\mathbb{G}_{\R^N}[Qu].$$

Now we  prove the claim (\ref{2.04}) by contradiction. We may assume that for some $x_0 \in \R^N\setminus\{0\}$,
\begin{equation}\label{5.1}
 \lim_{n\to+\infty }\int_{B_{n}(0)\setminus B_{\frac1{n}}(0)}\frac{Q(y)u(y)}{|x_0-y|^{N-2}}\,dy=+\infty.
\end{equation}

Let $\eta_0:[0,+\infty)\to[0,1]$ be a $C^\infty$ function such that $\eta_0$ equals to 1 in $[0,1]$ and
vanishes in $[2,+\infty)$. Since $Qu$ is locally H\"{o}lder continuous  in $\R^N\setminus\{0\}$, take
$$g_n(x)=\eta_0(\frac{|x|}n)(1-\eta_0(n|x|))Q(x)u(x),$$
then $g_n\le Qu$ and $g_n$ is H\"{o}lder continuous.  Let $w_n=\mathbb{G}_{B_{2n}(0)}[g_n]$, which  is the solution of
\begin{equation}\label{eq 3.1}
  \arraycolsep=1pt\left\{
\begin{array}{lll}
 \displaystyle    -\Delta u=g_n \quad
   &{\rm in}\quad B_{2n}(0),\\[2mm]
 \phantom{  -\Delta  }
 \displaystyle  u=0 &{\rm in}\quad \partial B_{2n}(0).
 \end{array}\right.
\end{equation}
We have that
$$\lim_{x\to0} w_n(x)|x|^{N-2}=0\quad{\rm and}\quad \liminf_{x\to0} u(x)|x|^{N-2}\ge0,$$
then by Lemma  \ref{lm cp} with $\mu=0$, it implies that
$$u\ge w_n\quad{\rm in} \quad B_{2n}(0)\setminus\{0\}.$$

We claim that there exists $c_3>0$ such that for any $n\ge 2$
\begin{equation}\label{2.101}
 G_{B_{2n}(0)}(x,y)\ge c_3|x-y|^{2-N},\quad\forall x,y\in B_{2n}(0),\, x\not=y.
\end{equation}
 In fact, it is known that
 $$ G_{B_2(0)}(x,y)\ge c_3|x-y|^{2-N},\quad\forall x,y\in B_1(0),\, x\not=y.$$
 It follows by the scaling property that for any $x,y\in B_n(0),\, x\not=y,$
 \begin{eqnarray*}
 G_{B_{2n}(0)}(x,y)  = n^{2-N} G_{B_2(0)}(x/n,y/n)
  \ge  c_3 |x-y|^{2-N}.
 \end{eqnarray*}

For $x_0\in \R^N\setminus\{0\}$, there exists $n_0>2$ such that $x_0\in B_{n_0}(0)\setminus B_{\frac1{n_0}}(0)$
thus, for $n\ge n_0$,
\begin{eqnarray*}
u(x_0) &\ge& w_n(x_0)= \int_{B_{2n}(0)}  G_{B_{2n}(0)}(x_0,y) g_n(y) dy \\
   &>& c_3\int_{B_{n}(0)\setminus B_{\frac1{n}}(0)}\frac{Q(y)u(y)}{|x_0-y|^{N-2}}dy
   \\&\to&+\infty\quad {\rm as}\quad n\to+\infty,
\end{eqnarray*}
which contradicts the fact that $u$ is a classical solution of (\ref{eq}). \hfill$\Box$

\bigskip

\noindent{\bf Proof of Proposition \ref{classification}.} Let $u_0$ be a classical solution of (\ref{eq 1.1}) such that $u_0\ge0$ in $B_r(0)\setminus\{0\}$, if $u_0$ is nontrivial, by Maximum Principle, we have that
$$u_0> 0\quad \quad {\rm in}\quad \R^N\setminus\{0\}.$$
Without loss of the generality, we put $r=1$.

We observe that $u_0$ is a solution of
$$\mathcal{L}_\mu u_0(x)=\mu (V(x)-|x|^{-2})u_0(x),\quad\forall\, x\in\R^N\setminus\{0\},$$
where $f_0:=\mu (V(x)-|x|^{-2})u_0(x)\ge0$ in $B_1(0)$

{\it Now we prove that $f_0\in L^1_{loc}(\R^N,\,d\mu)$.}  If not, we may assume that
$$\lim_{r\to0^+}\int_{B_1(0)\setminus B_r(0)}f_0 \, dx=+\infty.$$
By the fact $f_0\in C^\gamma(  B_1(0)\setminus \{0\})$, we have that
$$
 \lim_{r\to0^+} \int_{B_1(0)\setminus B_r(0)}f_0\, d\mu  =+\infty,
$$
then  there exists $R_n\in (0,1)$ such that
$$
  \int_{B_1(0)\setminus B_{R_n}(0)}f_0\, d\mu =n.
$$

\textbf{Let $w_n$ be the solution of
$$
 \arraycolsep=1pt\left\{
\begin{array}{lll}
 \displaystyle  \mathcal{L}_\mu u(x)= f_0(x)\eta_0(x)(1-\eta_0(2n|x|))\qquad
   &\forall x\in B_4(0)\setminus \{0\},\\[1.5mm]
 \phantom{  L_\mu   }
 \displaystyle  u(x)= 0\qquad  &\forall x\in \partial{ B_4(0)},\\[1.5mm]
 \phantom{   }
  \displaystyle \lim_{x\to0}u(x)\Phi_\mu^{-1}(x)=0,
 \end{array}\right.
$$
then it follows by comparison principle
$$u\ge w_n(x)\ge \int_{B_4(0)}\mathcal{G}_\mu(x,y) f_0(y)\eta_0(y)(1-\eta_0(2n|y|)) dy,$$
where $\mathcal{G}_\mu$ is the Green kernel of $\mathcal{L}_\mu $ in $B_4(0)\times B_4(0)$ subject to zero Dirichlet boundary condition.
It follows from \cite[Remark 4.1]{CQZ} that for $x,y\in B_3(0)\setminus\{0\}$, $x\not=y$,
$$\mathcal{G}_\mu(x,y)\ge
c[|x-y|^{2-N}+\frac{|x|^{\tau_+(\mu)} }{|x-y|^{N-2+\tau_+(\mu)}}+\frac{ |y|^{\tau_+(\mu)}}{|x-y|^{N-2+\tau_+(\mu)}}+ \frac{|x|^{\tau_+(\mu)}|y|^{\tau_+(\mu)}}{|x-y|^{N-2+2\tau_+(\mu)}}].$$
For $x_0\in \R^N$ with $|x_0|=2$ fixed,   it deduce that
\begin{eqnarray*}
u(x_0)\ge w_n(x_0) &\ge & c2^{\tau_+(\mu)}\int_{B_2(0)} \frac{|y|^{\tau_+(\mu)}}{|x_0-y|^{N-2+2\tau_+(\mu)}}f_0(y) dy  \\
  &\ge &  c\int_{B_1(0)\setminus B_{R_n}(0)}f_0\, d\mu\to+\infty\quad{\rm as}\quad n\to+\infty,
\end{eqnarray*}}
which is impossible. So we have that  $f_0\in L^1_{loc}(\R^N,\,d\mu)$.

From Lemma \ref{pr 2.1},  $u_0$ verifies that for some $k\ge0$,
$$
\int_{\R^N}u_0 \mathcal{L}_\mu^*(\xi)\, d\mu  = \int_{\R^N} \mu (V(x)-|x|^{-2})u_0(x)  \xi\, d\mu +c_\mu k\xi(0),\quad\forall\, \xi\in   C^{1.1}_0(\R^N),
$$
which implies (\ref{d 1.2}).

If $\mu\in(0,\mu_0)$, $f_0\ge0$ and $u_0\ge0$  in $\R^N\setminus\{0\}$, it deduces by (\ref{2.6})  that
$u_0\ge k\Phi_\mu.$
  \hfill$\Box$

\setcounter{equation}{0}
\section{Fundamental solutions}

\subsection{Existence}

To prove Theorem \ref{teo 2}, the following estimate plays an important role  in finding the fundamental solution.

\begin{lemma}\label{lm 3.1}
Assume that   $N\ge3$, $ 0< \mu<\mu_0$    and $V$ is a positive Hardy potential satisfying $(V)$ with $1\le \alpha<\frac{\mu_0}{\mu}$.
Let
$$\bar u_{k'}=\Phi_\mu+k'\Phi_{\mu'},$$
where $\mu'\in ( \mu\alpha,\,\mu_0)$ and $k'>0$.

Then   there exists $k^*>0$ depending on $\alpha,\, \mu$  such that $k'\ge k^*$,  $\bar u_{k'}$ is a super solution of
$$\mathcal{L}_{\mu V} u=0\quad{\rm in}\quad \R^N\setminus\{0\}.$$
\end{lemma}
{\bf Proof.}  From $(V)$, we have that $\frac{\mu'}{\mu}>\alpha$, so let $\alpha'\in(\alpha,\, \frac{\mu'}{\mu})$, there exists $r'>0$ such that
$$V(x)\le \arraycolsep=1pt\left\{
\begin{array}{lll}
\alpha' |x|^{-2}\quad&{\rm for}\quad |x|\ge r',\\[1mm]
 \phantom{   }
 |x|^{-2}+c_0\qquad&{\rm for}\quad 0<|x|<r'.
 \end{array}\right.
$$
There exists $\iota>0$ such that for $0<|x|<r'$,
$$\Phi_{\mu'}(x)\le \Phi_{\mu}(x)\le \iota |x|^{-2}\Phi_{\mu'}(x) .$$ Take  $r'=\min\{\sqrt{\frac{\mu'-\mu}{2c_0\mu}},1\}$, then
 we have that
 $$\frac{ \mu'- \mu}{|x|^2}-c_0 \mu\ge \frac{ \mu'- \mu}{2|x|^2}. $$
Let $k_0=\frac{ 2c_0\mu \iota }{ \mu'-\mu}$, and  for $k'\ge k_0$ and $x\in B_{r'}(0)\setminus\{0\}$, one implies that
\begin{eqnarray*}
-\Delta \bar u_{k'}(x)-\mu V\bar u_{k'}(x) &\ge & -c_0 \mu    \Phi_\mu(x) + [\frac{ \mu'- \mu}{|x|^2}-c_0 \mu  ]k \Phi_{\mu'}(x)
\\ &\ge& \left[ \frac{ \mu'- \mu}{2 }k'-c_0 \mu    \iota  \right] \frac{ \Phi_{\mu'}(x)}{|x|^2}
\\&\ge&0.
\end{eqnarray*}

 For $|x|\ge r'$, $V(x)\le \alpha' |x|^{-2}$ and there exists $\iota'>0$ such that
  $\Phi_{\mu'}\ge \iota' \Phi_{\mu}$. We note that $\iota=1$ if $r'=1$. Letting $k_1=\frac{(\alpha'-1)\mu  }{( \mu'-\alpha'\mu)\iota'}$, where  $\alpha' \mu-\mu'<0$.
  For $k'\ge k_1$ and $x\in \R^N\setminus B_{r'}(0)$
\begin{eqnarray*}
-\Delta \bar u_{k'}(x)-\mu V\bar u_{k'}(x)&\ge &-\frac{\mu(\alpha'-1) }{|x|^2}  \Phi_\mu(x) - \frac{(\alpha' \mu-\mu')k'}{|x|^2}  \Phi_{\mu'}(x)
\\ &\ge&  [(\mu'-\alpha' \mu)\iota' k' -(\alpha'-1)\mu   ]\frac{ \Phi_\mu(x)}{|x|^2}
\\&\ge&0.
\end{eqnarray*}
Therefore, taking $k^*=\max\{k_0,k_1\}$, the function $\bar u_{k'}$ with $k'\ge k^*$ is a super solution of
$$\mathcal{L}_{\mu V} \bar u_{k'}\ge0\quad{\rm in}\quad \R^N\setminus\{0\},$$
which ends the proof. \hfill$\Box$\medskip

\begin{lemma}\label{lm 2.1}
  Assume that  $N\ge3$, $k\ge0$, $ 0< \mu<\mu_0$, $V$ is a potential verifying $(V)$,  $\Omega$ is a $C^2$, bounded domain containing the origin and
  $\bar u_{k'}$ is the function constructed in Lemma \ref{lm 3.1}.
Then the problem
 \begin{equation}\label{eq 1.1f}
 \arraycolsep=1pt\left\{
\begin{array}{lll}
 \displaystyle   \mathcal{L}_{\mu V} u= 0\qquad
  & {\rm in}\quad  {\Omega}\setminus \{0\},\\[1.5mm]
    \phantom{    \mathcal{L}_{\mu V} }
    u=0& {\rm on}\quad  \partial{\Omega},
   \\[1.5mm]
 \displaystyle \lim_{x\to0}u(x)\Phi_\mu^{-1}(x)=k
 \end{array}\right.
\end{equation}
admits a unique solution  $v_k$, which satisfies
$$v_k\le k \bar u_{k'}\quad{\rm in}\quad \Omega\setminus\{0\}$$
and verifies  the $d\mu$-distributional identity
\begin{equation}\label{4.1.1}
\int_{\Omega} v_k \mathcal{L}^*_{\mu V}   \xi \,d\mu =  c_\mu\xi(0),  \quad \forall\, \xi\in C^\infty_c(\Omega).
\end{equation}
\end{lemma}
{\bf Proof.} From the linearity of $\mathcal{L}_{\mu V}$, we only have to prove this lemma with $k=1$.
Let $w_0$ be the solution of
$$ \arraycolsep=1pt\left\{
\begin{array}{lll}
 \displaystyle   \mathcal{L}_{\mu} u=c_\mu \delta_0\qquad
  & {\rm in}\quad   \Omega,\\[1.5mm]
    \phantom{    \mathcal{L}_{\mu } }
    u=0& {\rm on}\quad  \partial{\Omega}.
 \end{array}\right.$$
From Lemma \ref{lm cp}, we have that $0<w_0\le \bar u_{k'}$ in $\Omega\setminus\{0\}$.
From Theorem 1.3 in \cite{CQZ}, problem
$$ \arraycolsep=1pt\left\{
\begin{array}{lll}
 \displaystyle   \mathcal{L}_{\mu} u= \mu(V-|x|^{-2})w_0   \qquad
  & {\rm in}\quad   \Omega\setminus\{0\},\\[1.5mm]
    \phantom{    \mathcal{L}_{\mu } }
    u=0& {\rm on}\quad  \partial{\Omega},
   \\[1.5mm]
 \displaystyle \lim_{x\to0}u(x)\Phi_\mu^{-1}(x)=1
 \end{array}\right.$$
has a unique solution $w_1$. By Lemma \ref{lm cp}, we have that
$$w_0\le w_1\le \bar u_{k'}\quad{\rm in}\quad\Omega\setminus\{0\}.$$
Inductively,  for given $w_{n-2}\le w_{n-1}\le \bar u_{k'}$ in $\Omega$,  problem
$$ \arraycolsep=1pt\left\{
\begin{array}{lll}
 \displaystyle   \mathcal{L}_{\mu} u= \mu(V-|x|^{-2})w_{n-1}   \qquad
  & {\rm in}\quad   \Omega\setminus\{0\},\\[1.5mm]
    \phantom{    \mathcal{L}_{\mu } }
    u=0& {\rm on}\quad  \partial{\Omega},
   \\[1.5mm]
 \displaystyle \lim_{x\to0}u(x)\Phi_\mu^{-1}(x)=1
 \end{array}\right.$$
has a unique classical solution of $w_n$, which satisfies
$$w_{n-1}\le w_n\le \bar u_{k'}\quad{\rm in}\quad\Omega\setminus\{0\}$$
and by \cite[(1.12)]{CQZ},
\begin{equation}\label{4.2.3.1}
\int_{\Omega} w_n \mathcal{L}^*_{\mu }   \xi \,d\mu =\mu \int_{\Omega} [V(x)-|x|^{-2}] w_{n-1}\xi\, d\mu+ c_\mu\xi(0),  \quad \forall\, \xi\in C^\infty_0(\Omega).
\end{equation}
 Therefore, the sequence $\{w_n\}_n$ is convergent in $L^1(\Omega,\,d\mu)$. Let
$$v_{1}:=\lim_{n\to\infty} w_n,$$
then $$w_0 \le v_1\le \bar u_{k'}$$
and by the standard regularity result, it is known that $v_1$ is a classical solution of (\ref{eq 1.1f}), satisfying (\ref{4.1.1}) with $k=1$ by passing the limit of (\ref{4.2.3.1}) as $n\to+\infty$.
\smallskip\smallskip

\noindent{\bf Proof of Theorem \ref{teo 2}.}  Since $\mathcal{L}_{\mu V}$ is a linear operator, then we only have to prove the existence of Fundamental solution verifying (\ref{d 1.2}) with $k=1$. To this end,  denote by
$v_n$ the unique solution of (\ref{eq 1.1f}) with $k=1$ and $\Omega=B_n(0)$,  extend it by zero in $B_n^c(0)$ and still denote it by $v_n$.

{\it Claim 1: Letting $w_n=v_n-v_{n-1}$, then $w_n\ge 0$ on $\partial B_{n-1}(0)$. } Indeed,
 $$\arraycolsep=1pt\left\{
\begin{array}{lll}
 \displaystyle \mathcal{L}_{\mu V} w_n = 0\qquad
   {\rm in}\quad   B_{n-1}(0) \setminus \{0\},\\[1.5mm]
 \phantom{\mathcal{L}_{\mu V}    }
 \displaystyle  w_n\ge 0 \qquad  {\rm   on}\quad \partial B_{n-1}(0),\\[1.5mm]
 \phantom{   }
  \displaystyle \lim_{x\to0}w_n(x)\Phi_\mu^{-1}(x)=0,
 \end{array}\right.$$
  then for any $\epsilon>0$, there exists $r_\epsilon>0$ converging to zero as $\epsilon\to0$ such that
 $$w_n\ge -\epsilon v_{n-1}\quad{\rm in}\quad \partial B_{r_\epsilon}(0).$$
We see that
$$w_n\ge 0>-\epsilon \Phi_\mu \quad{\rm on}\quad \partial B_n(0),$$
then by Lemma \ref{lm ocp}, we have that
$$w_n\ge -\epsilon \Phi_\mu\quad{\rm in}\quad B_n(0)\setminus\{0\}. $$
By the arbitrary of $\epsilon$, we have that $w_n\ge 0$   in $B_{n-1}(0) \setminus\{0\}.$

So we have that
$$v_n\ge v_{n-1},\quad{\rm in}\quad B_{n-1}(0),$$
that is,  the sequence $\{v_n\}_n$ is  increasing with respect to $n$.
Similar to Claim 1, we have that that  $v_n\le \bar u_{k'}$ in $\R^N\setminus\{0\}$, then the sequence $\{v_n\}_n$ is convergent in $\R^N\setminus\{0\}$ and in $L^1_{loc}(\R^N,\, d\mu)$. Let
$$u_{1}:=\lim_{n\to\infty} v_n,$$
then $$v_0 \le u_1\le \bar u_{k'}$$
and by the standard regularity result, it is known that $u_1$ is a classical solution of (\ref{eq 1.1}), satisfying (\ref{1.3}) with $k=1$ and (\ref{1.3.1}).

 For any $\xi\in C^\infty_c(\R^N)$, there exists $n_0$ such that
$\xi\in C^\infty_c(B_{n_0}(0))$, then from   (\ref{4.1.1}),  we have that for $n\ge n_0$
$$
\int_{\R^N} v_n \mathcal{L}^*_{\mu V}   \xi \,d\mu  =  c_\mu\xi(0),
$$
passing to the limit as $n\to+\infty$, then
 $u_1$ verifies  (\ref{d 1.2}) with $k=1$.

Finally, we   prove that $u_{1}$ is the minimal solution of (\ref{eq 1.1}) verifying (\ref{d 1.2}) with $k=1$.   Indeed,  let $u>0$ be a fundamental solution of (\ref{eq 1.1}) verifying (\ref{d 1.2}) with $k=1$, then  from Lemma \ref{lm ocp}, we have that
\[ u_1  \ge v_n,\]
which implies $u\ge u_1$ by the fact that $u_1=\lim_{n\to+\infty} v_n$.
\hfill$\Box$



\subsection{Nonexistence}

We prove the nonexistence of fundamental solutions of (\ref{eq 1.1}) by contradiction. Assume that $u$ is a positive fundamental solution to problem (\ref{eq 1.1}) and we will obtain a contradiction from its decay at infinity.

The following observation is very important for proving the nonexistence of  fundamental  solutions of $\mathcal{L}_\mu$.

\begin{lemma}\label{lm 2.2}
For any $e\in \mathcal{S}^{N-1}$ and $\mu\in(\mu_0,0)$, we have that
 \begin{equation}\label{2.2}
\mu c_N\int_{\R^N} \frac{\Phi_\mu(y)}{|e-y|^{N-2}|y|^2}\,dy =1.
 \end{equation}
\end{lemma}
{\bf Proof.} We have that $\Phi_{\mu}$ is a classical solution of (\ref{eq 1.2})
and, by   Lemma \ref{re 2.1}, it infers that $\Phi_{\mu}$ could be expressed by
\begin{eqnarray*}
 \Phi_\mu(x) =\mathbb{G}_{\R^N}[\mu |\cdot|^{-2}\Phi_\mu](x)
    =  c_N \mu  \int_{\R^N} \frac{\Phi_\mu(y)|y|^2}{|x-y|^{N-2}}\,dy,
\end{eqnarray*}
which implies that
\begin{equation}\label{eq 4.1}
 \Phi_\mu=\mu \mathbb{G}_{\R^N}[|\cdot|^{-2}\Phi_\mu]\quad{\rm in}\quad \R^N\setminus\{0\}
\end{equation}
and (\ref{2.2}) by taking $x\in \partial B_1(0)$.\hfill$\Box$\medskip

We remark that  from the view of \cite{BDT}, the function $\Phi_\mu$ is a distributional solution of
\begin{equation}\label{eq 2.1}
 -\Delta u=\frac{\mu}{|x|^2} u\qquad{\rm in}\quad \mathcal{D}'(\R^N).
\end{equation}

We note that the fundamental solution $u$ of (\ref{eq 1.1})  verifies (\ref{d 1.2}) with some $k>0$, then from the proof of Theorem \ref{teo 2}, there exists a minimal fundamental solution
 $u_k$ of (\ref{eq 1.1})  verifying (\ref{d 1.2}) with such $k>0$ and from Lemma \ref{re 2.1}, it has the formula
\begin{equation}\label{4.4}
 u_k= \mu \mathbb{G}_{\R^N}[Vu_k].
\end{equation}

\noindent{\bf Proof of Theorem \ref{teo 1}.} We prove Theorem \ref{teo 1} by contradiction. Assume that $u$ is a fundamental solution of (\ref{eq 1.1}) verifying (\ref{d 1.2}) with some $k>0$ and $u_k$ is the minimal fundamental solution.
Since $u_k>0$ and by Assumption $(\widetilde{V})$, we have that
$$f(x):=\mu(V(x)-|x|^{-2}) u_k(x)\ge0,$$
 then by (\ref{2.6}),
 \begin{equation}\label{2.1}
u_k(x)\ge k\Phi_\mu(x),\quad\forall x\in \R^N\setminus\{0\}.
\end{equation}

Let $$\theta=\min\{(\frac{1+\beta}{2})^{\frac{1}{ N-2+\tau_-(\mu)}},1-2^{\frac1{2-N}}\}$$ and for any $\epsilon\in(0,\frac{\beta-1}{\beta+1})$, there exists $r_0>1$ such that
$$V(x)|x|^2\ge \beta(1-\epsilon)\quad{\rm for}\quad   |x|\ge r_0.$$
Taking $r_1=\theta^{-1}r_0$, it implies by (\ref{eq 2.1}) that for $x\in B_{r_1}^c(0)$,
\begin{eqnarray*}
u_k(x) &\ge& k \mu  \beta(1-\epsilon) c_N \int_{\R^N\setminus B_{r_0}(0) } \frac{\Phi_\mu(y)|y|^{-2}}{|x-y|^{N-2}}dy \\
    &\ge&k \mu  \beta(1-\epsilon) c_N  \int_{\R^N\setminus B_{\frac{r_0}{r_1}}(0) } \frac{\Phi_\mu(y)|y|^{-2}}{|e_x-y|^{N-2}}dy\, \Phi_\mu(x) \\
    &\ge& k \sigma \Phi_\mu(x),
\end{eqnarray*}
where  $|e_x-y|^{2-N}\le (1-\theta)^{2-N}\le 2$ and
\begin{eqnarray*}
\sigma &=&  c_N\beta \mu (1-\epsilon) \left[\frac{1}{c_N\mu}- \int_{ B_{\theta}(0) } \frac{\Phi_\mu(y)|y|^{-2}}{|e_x-y|^{N-2}}dy\right] \\
   &\ge&[\beta-2\theta^{ N-2+\tau_-(\mu)}] (1-\epsilon)
      \\&\ge&\frac{1+\beta}{2}(1-\epsilon)>1
\end{eqnarray*}
by the choice of $\epsilon$.

 Denote $r_n=\theta^{-n}r_0$, then we may assume that
 $$u_k(x)\ge k \sigma^{n-1}  \Phi_\mu(x), \quad\forall x\in B_{r_{n-1}}^c(0).$$
We observe that   for $x\in B_{r_{n}}^c(0)$,
\begin{eqnarray*}
u_k(x) &\ge&  k\mu c_N(1-\epsilon) \sigma^{n-1}\int_{\R^N\setminus B_{r_{n-1}}(0) } \frac{ \Phi_\mu(y)|y|^{-2}}{|x-y|^{N-2}}dy \\
    &\ge& k \mu  c_N(1-\epsilon) \sigma^{n-1} \int_{\R^N\setminus B_{\theta}(0) } \frac{ \Phi_\mu(y) |y|^{-2}}{|e_x-y|^{N-2}}dy\,  \Phi_\mu(x)\\
    &\ge& k \sigma^{n}  \Phi_\mu(x).
\end{eqnarray*}
As a conclusion, we have that
\begin{equation}\label{2.3}
  u_k(x)\ge k \sigma^{n} \Phi_\mu(x), \quad\forall x\in B_{r_{n}}^c(0).
\end{equation}

Finally, we shall get a contradiction by the decay at infinity. We claim that there exists $\mu^*>\mu_0$ such that
\begin{equation}\label{2.4}
\sigma\theta^{N-2+\tau_-(\mu)}\ge1.
\end{equation}
Indeed,
\begin{eqnarray*}
\sigma\theta^{N-2+\tau_-(\mu)}   &\ge&   \left[\min\{\beta\frac{1+\beta}{2},\, \beta(1-2^{\frac1{2-N}})\}-2\right](1-\epsilon)
\end{eqnarray*}
So there exists $\beta^*>1$ such that for $\beta>\beta^*$, one has that
$$\left[\min\{\beta\frac{1+\beta}{2},\, \beta(1-2^{\frac1{2-N}})\}-2\right]>1$$
and then we may choose $\epsilon>0$ small enough, we have that
\begin{equation}\label{2.5}
\sigma\theta^{N-2+\tau_-(\mu)}\ge1.
\end{equation}
Then fix some point $x\in \partial B_1(0)$ and    there exists $c_4>0$ such that for $y\in B_{r_0}^c(0),$
$$|x-y|^{2-N}\ge c_4|y|^{2-N}$$
and then we have that
\begin{eqnarray*}
u_k(x) &\ge&   \mu  c_N (1-\epsilon)\int_{\R^N\setminus B_{r_0}(0) } \frac{u(y)|y|^{-2}}{|x-y|^{N-2}}dy \\
    &\ge&  k\mu  c_N  (1-\epsilon) \sum_{n=2}^{+\infty}\left(\sigma^n\int_{B_{r_n}(0)\setminus B_{r_{n-1}}(0) } \frac{ \Phi_\mu(y) |y|^{-2}}{|x-y|^{N-2}}dy \right) \\
    &\ge&  c_4k\mu  c_N(1-\epsilon) \sum_{n=2}^{+\infty} \left(\sigma^n\int_{B_{r_n}(0)\setminus B_{r_{n-1}}(0) } \frac{ \Phi_\mu(y) |y|^{-2}}{|y|^{N-2}}dy \right)
    \\&=& c_4 k\frac{N-2}{2} (1-\epsilon)(1-\theta^{ \frac{N-2}2} )  r_0^{\frac{2-N}2} \sum_{n=2}^{+\infty}\left(\sigma\theta^{ N-2+\tau_-(\mu) }\right)^n
    \\&=&+\infty,
\end{eqnarray*}
which is impossible.
The proof ends.\hfill$\Box$

\subsection{Proof of Corollary \ref{cr 1.1} }

\noindent{\bf Proof of Corollary \ref{cr 1.1}.} From (\ref{Vr}), we know that
$$0\le V_\rho(x)-|x|^{-2}\le \rho-1,\quad \forall\,x\in \R^N\setminus\{0\}$$
and
$$\lim_{|x|\to+\infty} V_\rho(x)|x|^2=\rho.$$
For some $\rho^*$, $V_{\mu\rho}$ verifies $(\tilde V)$  if $\rho> \rho^*$, then it deduces by  Theorem \ref{teo 1} that $\mathcal{L}_{\mu V_\rho}$ has no positive fundamental solutions.

So $V_\rho$ verifies the assumption $(V)$ if $\rho<\frac{\mu_0}{\mu}$,
then it deduces from Theorem \ref{teo 2}, there are positive fundamental solutions
for $\mathcal{L}_{\mu V_\rho}$.

Obviously, $\rho^* \ge \frac{\mu_0}{\mu}$. To finish the proof of Corollary \ref{cr 1.1}, we only prove the following argument:
Let $\rho_1>\rho_2\ge 1$, and if $\mathcal{L}_{\mu V_{\rho_1}}$ has positive fundamental solutions, then
 $\mathcal{L}_{\mu V_{\rho_2}}$ has fundamental solutions.

Indeed,  let $u_{\rho_1}$ be a fundamental solution of (\ref{eq 1.1}) verifying (\ref{d 1.2}) with $k=1$, replacing $V$ by $V_{\rho_1}$.

  We  recall that
$v_n$ is the unique solution of (\ref{eq 1.1f}) with $k=1$ and $\Omega=B_n(0)$,  extend it by zero in $B_n^c(0)$ and still denote it by $v_n$.

Since $V_{\rho_2}\ge |x|^{-2}$,  $\{v_n\}_n$ is an increasing sequence and   $V_{\rho_2}<V_{\rho_1}$ implies that the function  $u_{\rho_1}$ is an upper bound for  $\{v_n\}_n$.
As a conclusion,  the limit of $\{v_n\}_n$ as $n\to+\infty$ is a positive fundamental solution $\mathcal{L}_{\mu V_{\rho_2}}$. The other details could see the proof of  Theorem \ref{teo 2}.
\hfill$\Box$

 \setcounter{equation}{0}
 \section{  Discussion  }

\subsection{ On $0<\mu<\mu_0$ and $V(x)\le |x|^{-2}$ }

We discuss the existence of the fundamental solutions for $\mathcal{L}_{\mu V}$ when $V(x)\le |x|^{-2}$.
\begin{theorem}\label{teo 5.1}
Assume that   $N\ge3$, $ \mu\in(0,\mu_0)$   and  $V\in C^\gamma_{loc}(\R^N\setminus\{0\})$ with $\gamma\in(0,1)$ is a positive potential satisfying that for some $c_5>0$,
\begin{equation}\label{v5.1}
 |x|^{-2} -c_5\le V(x)\le |x|^{-2},\quad \forall x\in \R^N\setminus \{0\}.
\end{equation}
Then  $\mathcal{L}_{\mu V}$
 has positive  fundamental solutions and for  $k>0$, there is a  minimal positive solution $u_k$ verifying (\ref{d 1.2}) with such $k$ and satisfying  that
 \begin{equation}\label{6.1}
\lim_{x\to0}u_k(x)\Phi_\mu^{-1}(x)=k.
 \end{equation}

$(i)$ If there is some $\nu\in(0,1)$,  $$\lim_{x\to+\infty} V(x)|x|^{2}=\nu,$$
then the minimal fundamental solution $u_k$ verifies that for any $\mu'\in (\nu\mu,\, \mu)$, there exists  $c_6>0$ such that
  \begin{equation}\label{6.2}
 u_k\le c_6\Phi_{\mu'}\quad {\rm in}\quad  B_1^c(0).
 \end{equation}

$(ii)$ If  there exists $\tau>2$ such that
\begin{equation}\label{6.4}
 \limsup_{x\to+\infty} V(x)|x|^{\tau}<+\infty,
\end{equation}
then the minimal fundamental solution $u_k$ verifies that   there exists  $c_7>0$ such that
  \begin{equation}\label{6.3}
 u_k\le c_7\Phi_{0}\quad {\rm in}\quad   B_1^c(0).
 \end{equation}

 \end{theorem}

Before proving Theorem \ref{teo 5.1}, we have to use the following lemma.

\begin{lemma}\label{lm 5.1}
Assume that  $N\ge3$, $ \mu\in(0,\mu_0]$,  $V\in C^\gamma_{loc}(\R^N\setminus\{0\})$ with $\gamma\in(0,1)$ is a positive potential satisfying  (\ref{v5.1}) and
 $\Omega$ is a $C^2$, bounded domain containing the origin.
Then  the problem
\begin{equation}\label{eq 6.1}
 \arraycolsep=1pt\left\{
\begin{array}{lll}
 \displaystyle   \mathcal{L}_{\mu V} u =0\qquad
   &{\rm in}\quad  {\Omega\setminus\{0\}},\\[1.5mm]
 \phantom{   \mathcal{L}_{\mu V}  }
 \displaystyle  u= 0\qquad  &{\rm   on}\quad \partial{\Omega},\\[1.5mm]
 \phantom{    }
 \displaystyle\lim_{x\to0}u(x)\Phi_\mu^{-1}=1
 \end{array}\right.
\end{equation}
has a unique positive solution $v$ satisfying that
\begin{equation}\label{6.7}
  \int_{ \Omega} u \mathcal{L}_{\mu V}^* \xi  d\mu   =c_\mu k\xi(0),\quad\forall\, \xi\in C^{1.1}_0(\Omega).
\end{equation}

\end{lemma}
{\bf Proof.}
{\it The existence.} We observe that $\mathcal{G}_\mu$ is the  solution of  $\mathcal{L}_\mu u=c_\mu\delta_0$ in $\Omega$ in the $d\mu$-distributional sense, subjecting to $u=0$ on $\partial \Omega$. Then
$$\mathcal{G}_\mu=\mu\mathbb{G}_{\Omega}[|\cdot|^{-2}\mathcal{G}_\mu ],$$
where $\mathbb{G}_{\Omega}$ is the Green's operator defined by the Green kernel $G_{\Omega}$ of $-\Delta$ in $\Omega\times \Omega$.

The existence of solution of (\ref{eq 6.1}) could be approximated by the sequence
$$w_0=\mathcal{G}_\mu,\quad{\rm and}\quad w_n=\mu\mathbb{G}_{\Omega}[Vw_{n-1}].$$
It follows by
$$w_1=\mu\mathbb{G}_{\Omega}[Vw_0]\le \mu\mathbb{G}_{\Omega}[|\cdot|^{-2}w_0]=w_0,$$
so inductively, we have that $\{w_n\}_n$ is a decreasing sequence. Furthermore, $\{w_n\}_n$ is a positive sequence.

When $\mu<\mu_0$, let $t=\max\{1,\, \frac{c_0\mu}{\mu_0-\mu}\}$ and denote
$$w_t(x)=\Phi_\mu(x)-t|x|^{-\frac{N-2}{2}}\quad{\rm for} \ \ 0<|x|<r_t:=t^{\frac1{\tau_-(\mu)+\frac{N-2}{2}}},$$
then
\begin{eqnarray*}
\mathcal{L}_{\mu V}w_t(x) &=&  \mathcal{L}_\mu w_t(x)+\mu (|x|^{-2}-V(x))w_t(x)\\
   &\le&  -t(\mu_0-\mu) |x|^{-\frac{N-2}{2}-2}+c_0\mu|x|^{\tau_-(\mu)}
   \\&\le &0.
\end{eqnarray*}

When $\mu=\mu_0$, let $t=\max\{2,\,  8c_0\mu_0 \}$ and denote
$$w_t(x)=\Phi_{\mu_0}(x)-t|x|^{-\frac{N-2}{2}}(-\log|x|)^{\frac12}\quad{\rm for} \ \ 0<|x|<\frac14,$$
then
\begin{eqnarray*}
\mathcal{L}_{\mu_0 V}w_t(x) &=&  \mathcal{L}_{\mu_0} w_t(x)+\mu_0 (|x|^{-2}-V(x))w_t(x)\\
   &\le&  -\frac14t  |x|^{-\frac{N-2}{2}-2}(-\log |x|)^{-\frac12}+c_0\mu_0|x|^{-\frac{N-2}2}(-\log |x|)
   \\&\le &0.
\end{eqnarray*}

Since $w_0>0$ on $\partial B_{r_t}(0)$, then from Lemma \ref{lm cp}, we have that
$$w_1\ge w_t\quad{\rm in}\quad   B_{r_t}(0),$$
which, inductively, implies that for any $n\ge 1$
$$w_n\ge w_t\quad{\rm in}\quad   B_{r_t}(0).$$

Thus, $\{w_n\}_n$ is convergent,  letting $v_k=\lim_{n\to+\infty} w_n$, we have that
$$v_k=\mathbb{G}_{\Omega}[Vv_k]$$
and
\begin{equation}\label{6.6}
 (\Phi_\mu-t\Gamma_{\mu_0})_+\le v_k\le \mathcal{G}_\mu,
\end{equation}
then
$$\mathcal{L}_{\mu V} v_k=0\quad{\rm in}\quad \Omega\setminus\{0\},$$
Integrate over $\Omega\setminus B_r(0)$ and pass the limit as $r\to0^+$, then we deduces that $v_k$ is the $d\mu$-distributional solution of
(\ref{eq 6.1}) by using (\ref{6.6}), which also implies (\ref{6.7}). The calculations could refer to the proof of Theorem 1.1 in \cite{CQZ}.

{\it The uniqueness.}   Let $u_1,u_2$ two solutions of (\ref{eq 6.1}), let $w=u_1-u_2$ and then
 $w$ verifies that for any $\xi\in C^{1.1}_0(\Omega)$, $\xi\ge0$,
\begin{equation}\label{7.1}
 \int_{\Omega} |w|\mathcal{L}_\mu^*\xi\, d\mu +\mu\int_{\Omega}(|x|^{-2}-V) {\rm sign}(w) w\xi\, d\mu\le 0,
\end{equation}
thanks to the fact that
$$\mathcal{L}_{\mu V}=\mathcal{L}_{\mu}+\mu (|x|^{-2}-V),$$
where $\mu (|x|^{-2}-V)\ge0$.
Taking $\xi$ the solution of  $\mathcal{L}_\mu^*u=1$ in $\Omega$, subjecting to $u=0$ on $\partial \Omega$,
we derive that
$$
 \int_{\Omega} |w|\, d\mu\le 0,
$$
which implies that $w=0$ a.\,e. in $\Omega$, and the uniqueness follows.\smallskip
\hfill$\Box$

\medskip

\noindent{ \bf Proof of Theorem \ref{teo 5.1}.} Let $v_n$ be the unique solution of (\ref{eq 6.1}) with $\Omega=B_n(0)$,  extending $v_n$ by zero in $\R^N\setminus \overline{B_n(0)}$, still denoting by $v_n$. Let $w_n=v_n-v_{n-1}$, then $w_n\ge 0$ on $\partial B_{n-1}(0)$ by the observation  that
$$ v_{n}\ge 0= v_{n-1}\quad{\rm on}\quad \partial B_{n-1}(0).$$

{\it We claim that for any $n\ge 2$
$$v_n\ge v_{n-1}\quad{\rm in}\quad B_{n-1}(0),$$
that is,  the sequence $\{v_n\}_n$ is  increasing with respect to $n$.}

In fact, let $w=v_{n-1} -v_n$ be a solution of
 $$\arraycolsep=1pt\left\{
\begin{array}{lll}
 \displaystyle \mathcal{L}_{\mu V} u \le 0\qquad
   {\rm in}\quad  B_{n-1}(0)\setminus \{0\},\\[1.5mm]
 \phantom{ \mathcal{L}_{\mu V}    }
 \displaystyle  u\le 0 \qquad  {\rm   on}\quad \partial B_{n-1}(0),\\[1.5mm]
 \phantom{   }
  \displaystyle \lim_{x\to0}u(x)\Phi_\mu^{-1}(x)=0,
 \end{array}\right.$$
  then, together with
  $$\lim_{x\to0} v_{n-1}(x) \Phi_\mu^{-1}(x)=1,$$ for any $\epsilon>0$, there exists $r_\epsilon>0$ converging to zero as $\epsilon\to0$ such that
 $$w\le \epsilon\, v_{n-1} \quad{\rm in}\quad \partial B_{r_\epsilon}(0).$$
We see that
$$w\le 0=\epsilon\, v_{n-1} \quad{\rm on}\quad \partial B_{n-1}(0),$$
then by Lemma 2.1 in \cite{CQZ}, we have that
$$w\le \epsilon\, v_{n-1} \quad{\rm in}\quad B_{n-1}(0)\setminus\{0\}. $$
By the arbitrary of $\epsilon$, we have that $w\le 0$   in $\Omega\setminus\{0\}.$ We complete the proof of the claim.\smallskip

We see that the sequence $\{v_n\}_n$ is convergent, since the upper bound  is  $\Phi_\mu$, from the fact $V(x)\le |x|^{-2}$ and $\mu<\mu_0$,
i.e.  $v_n\le \Phi_\mu$ in $\R^N\setminus\{0\}$, then the sequence $\{v_n\}_n$ is convergent in $\R^N\setminus\{0\}$ and in $L^1_{loc}(\R^N,\, d\mu)$. Let
$$u_{1}:=\lim_{n\to\infty} v_n,$$
then $$v_0 \le u_1\le \Phi_\mu\quad{\rm in}\quad \R^N\setminus\{0\}$$
and by the standard regularity result, it is known that $u_1$ is a classical solution of (\ref{eq 1.1}), satisfying (\ref{1.3}) with $k=1$.

 The function $u_k:=\lim_{n\to+\infty}v_n$ is a minimal positive fundamental solution verifying (\ref{d 1.2}) with such $k$.
See the proof of Theorem \ref{teo 2}.

To prove (\ref{6.2}), we have to find some suitable upper bound for the sequence $\{v_n\}$.
Since
$$\lim_{x\to+\infty} V(x)|x|^{2}=\nu,$$
then for any $\nu'>\nu$, there exists $r_{\nu'}$ such that
 $$V(x)\le \nu'|x|^{-2} \quad{\rm for}\quad |x|>r_{\nu'}$$
So $\Phi_{\mu'}$ with $\mu'=\nu'\mu>\nu\mu$ is a super solution of
$$-\Delta u=\mu V u\quad{\rm in}\quad \R^N\setminus B_{r_\nu'}(0).$$
Since $\{v_n\}$ is controlled by $k\Phi_\mu$, so there exists $c_8>0$ such that for any $n$
$$v_n\le c_8\Phi_{\mu'}\quad{\rm on}\quad \partial B_{r_\nu'}(0). $$
By Comparison Principle, we have that
$$u_k\le c_8\Phi_{\mu'}.$$

To prove (\ref{6.3}), the upper bound could be constructed by  $c_9(|x|^{2-N}-|x|^{4-N-\tau})$ in $  B_r^c(0)$, where
$c_9>0$, $r>1$ is such that
$$v_n\le c_9(r^{2-N}-r^{4-N-\tau})\quad{\rm on}\quad \partial B_{r}(0). $$
By choosing $c_9$ again, we obtain (\ref{6.3}).\hfill$\Box$

\begin{remark}
The authors in \cite{MP} shows that when $p\in(\frac{N}{N-2}, \frac{N+2}{N-2})$, the elliptic problem
$$-\Delta u=u^p\quad{\rm in}\quad \R^N\setminus \{0\},$$
has a sequence of the fast decay solutions  $\{u_k\}_k$
such that
$$u_k(x)\sim \arraycolsep=1pt\left\{
\begin{array}{lll}
 \displaystyle c_p |x|^{-\frac{2}{p-1}}\qquad
  & {\rm at\ the\ origin},\\[1mm]
 \phantom{  }
 \displaystyle  k|x|^{2-N}\qquad   &{\rm   at\ infinity},
 \end{array}\right. $$
 where
 $$c_p=[\frac{2}{p-1}(N-2-\frac{2}{p-1})]^{\frac1{p-1}}.$$
 More related isolated singularities could refer to \cite{GS,P1}.

We observe that $0<c_p^{p-1}<\mu_0$ for $p \in(\frac{N}{N-2}, \frac{N+2}{N-2})$
 and $c_p^{p-1}=\mu_0$ if $p=\frac{N+2}{N-2}$.
So $u_k^{p-1}$ plays an  role of Hardy potential $\mu V$ in Theorem \ref{teo 5.1}.

\end{remark}

\subsection{ On  $\mu=\mu_0$}

In this subsection, we discuss   the fundamental solutions for $\mathcal{L}_{\mu V}$ when $\mu=\mu_0$. Precisely, we have the following result.

\begin{theorem}\label{teo 5.2}

Assume that   $N\ge3$, $ \mu=\mu_0$   and positive potential  $V\in C^\gamma_{loc}(\R^N\setminus\{0\})$ with $\gamma\in(0,1)$ is a positive potential satisfying that
 $$ |x|^{-2}-c\le  V(x)\le |x|^{-2} ,\quad \forall\, x\in\R^N\setminus\{0\}\quad{\rm and}\quad
  \limsup_{x\to0} V(x)|x|^2=\varrho<1,$$
then $\mathcal{L}_{\mu_0 V}$ has  positive fundamental solutions and for  $k>0$, there is a  minimal positive solution $u_k$ verifying (\ref{d 1.2}) with such $k$. Furthermore,
 \begin{equation}\label{6.5}
\lim_{x\to0}u_k(x)\Phi_{\mu_0}^{-1}(x)=k
 \end{equation}
 and for any $\mu'\in(\varrho \mu_0, \mu_0)$, there exists $c_{10}>0$ such that
 $$  u(x)\le c_{10}\Phi_{\mu'}(x),\quad \forall\, x\in  B_1^c(0).$$

 \end{theorem}
{\bf Proof.} From Lemma \ref{lm 5.1},   Let $v_n$ be the unique solution of (\ref{eq 6.1}) with $\Omega=B_n(0)$,  extending $v_n$ by zero in $\R^N\setminus \overline{B_n(0)}$, still denoting by $v_n$. Similar to the proof of the Theorem \ref{teo 5.1}, we have that $\{v_n\}$ is an increasing sequence.

 We only have to construct a super bound for this sequence.
Since $\varrho<1$ and $\varrho \mu_0<\mu_0$, then for any $\varrho'\in(\varrho,\, 1)$, there exists $r'\ge2$ such that $$V(x)\le \varrho'|x|^{-2}\quad{\rm for} \quad |x|>r',$$
so for $\varrho''\in(\varrho',1)$,
$$-\Delta \Gamma_{\varrho'' \mu_0}-\mu_0 V \Gamma_{\varrho'' \mu_0}\ge \mu_0(\varrho''-\varrho')|x|^{-2} \Gamma_{\varrho'' \mu_0} \quad{\rm in}\quad     B_{r'}^c(0).$$
So we take
$$\Lambda(x)=[1-\eta_0(\frac{x}{r'}) ]\Gamma_{\varrho'' \mu_0}(x),\quad\forall\, x\in\R^N\setminus\{0\},$$
where $\eta_0:[0,+\infty)\to[0,1]$  is a decreasing, smooth function such that
$\eta_0(t)=0$ for $t>2$ and $\eta_0(t)=1$  for $0\le t\le1$.

There exists $K_1>0$ such that
$$\Phi_{\mu_0}+K_1\Gamma_{\mu_0}>0\quad {\rm in}\quad B_{2r'}(0).$$
Let
$$\bar u_{k',k''}=\Phi_{\mu_0}+k'\Gamma_{\mu_0}+k''\Lambda.$$

For $0<|x|\le r'$,
\begin{eqnarray*}
-\Delta \bar u_{k',k''}-\mu_0 V\bar u_{k',k''} &\ge & \mu_0[|x|^{-2}-V(x)]\Phi_{\mu_0}(x) +\mu_0[|x|^{-2}-V(x)]k' \Gamma_{ \mu_0}(x)
\\&\ge&0.
\end{eqnarray*}

For $|x|>2r'$, there exists $K_2\ge K_1$ such that for $k''\ge K_2$,
\begin{eqnarray*}
-\Delta \bar u_{k',k''}-\mu_0 V\bar u_{k',k''} &\ge & \mu_0[|x|^{-2}-V(x)]\Phi_{\mu_0}(x) +\mu_0[|x|^{-2}-V(x)]k' \Gamma_{ \mu_0}(x)\\&&+ \mu_0[\varrho''|x|^{-2}-V(x)]k''\Gamma_{\varrho''\mu_0}(x)
\\ &\ge& \frac{\mu_0 }{|x|^2}\left[(1-\varrho')  \Phi_{\mu_0}(x) +(\varrho'' -\varrho')k'' \Gamma_{\varrho''\mu_0}(x) \right]
\\&\ge&0.
\end{eqnarray*}

Now we fix $k''=K_2$. For $r'\le |x|\le 2r'$, since $-\Delta \Lambda-\mu_0 V\Lambda\ge c$ is bounded,   then  for $k'>K_1$ big enough
\begin{eqnarray*}
-\Delta \bar u_{k',k''}-\mu_0 V\bar u_{k',k''} &\ge & \mu_0[|x|^{-2}-V(x)]\Phi_{\mu_0}(x) +\mu_0[|x|^{-2}-V(x)]k' \Gamma_{ \mu_0}(x) - \mu_0k''c
\\ &\ge& \frac{\mu_0 }{|x|^2}(1-\varrho') (k'-K_1) \Gamma_{\mu_0}(x) -   \mu_0k''c
\\&\ge&0.
\end{eqnarray*}
Therefore, for suitable $k',k''$,  the function $\bar u_{k',k''}$ is a super solution of
$$\mathcal{L}_{\mu_0 V} \bar u_{k',k''}\ge0.$$
The proof ends.\hfill$\Box$

\setcounter{equation}{0}
 \section*{  Appendix }
 In this appendix, we prove the classification of the isolated singularities for  positive solutions of (\ref{eq 2.1.1}), which is motivated by Proposition 4.2 in \cite{CQZ} and
 paper \cite{BL}. \medskip

\noindent{\bf Proof of Lemma \ref{pr 2.1}.}  Let
$$ \bar u(r)=|\mathcal{S}^N|^{-1}r^{1-N}\int_{\partial B_r(0)} u(x) d\omega(x).$$
For $r\in(0,1)$, we have that
$$-\bar u''(r)-\frac{N-1}{r}\bar u'(r)-\frac{\mu}{r^2} \bar u(r)\ge \bar f(r),$$
where
$$\bar f(r)=r^{1-N}\int_{\partial B_r(0)} f(x) d\omega(x).$$
Denote $$v(r)=r^{-\tau_+(\mu)} \bar u(r),$$
then
$$-v''(r)-\frac{N+2\tau_+(\mu)-1}{r} v'(r)\ge r^{-\tau_+(\mu)} \bar f,$$
where $N+2\tau_+(\mu)$ plays the dimensional role. From $f\in L^1_{loc}(\R^N,\, d\mu)$, we have that
$$\int_0^{1}  r^{-\tau_+(\mu)} \bar f(r) r^{N+2\tau_+(\mu)-1}dr= \int_{B_{1}(0)}  |x|^{\tau_+(\mu)}|f(x)| dx  <+\infty. $$

From the step 1 in the proof of Theorem 1.1 in \cite{BL}, there exists $c_{11}>0$ such that
$$v(r)\le  \left\{\arraycolsep=1pt
\begin{array}{lll}
c_{11} |x|^{2-(N+2\tau_+(\mu))}\quad
   &{\rm if}\quad N+2\tau_+(\mu)\ge 3,\\[1.5mm]
 \phantom{   }
  c_{11}(-\ln|x|) \quad  &{\rm   if}\quad N+2\tau_+(\mu)=2,
 \end{array}
 \right.$$
that is, $u\le c_{11}\Phi_\mu$.
So for $\xi\in  C^\infty_c(\R^N)$, it is well-defined that
$$\int_{\R^N} u \mathcal{L}^*_\mu(\xi)\, d\mu <+\infty.$$
 We observe that for $\xi\in  C^\infty_c({\R^N}\setminus\{0\})$, it follows by Divergence theorem that
$$\int_{{\R^N}}u_{\mu}  \mathcal{L}^*_{\mu}(\xi)   \Gamma_{\mu}dx -\int_{{\R^N}}f\,\xi \Gamma_{\mu} dx=0.$$
By Schwartz Theorem (\cite[Theorem XXXV]{S}),  there exists a multiple index $p$,
$$u_{\mu}\Gamma_{\mu} \mathcal{L}^*_{\mu}-f\Gamma_{\mu}=\sum_{|a|=0}^p k_a D^{a}\delta_0,$$
i.e. for any $\xi\in C^\infty_c({\R^N})$
\begin{equation}\label{4.2}
 \int_{{\R^N}}u_{\mu} (\mathcal{L}^*_{\mu}\xi-f\xi)  \, d\mu =\sum_{|a|=0}^p k_a D^{a}\xi(0).
\end{equation}
We are left to show that $k_a=0$ for $|a|\ge 1$.
For multiple index $\bar a\not=0$, taking $\xi_{\bar a}(x)=x_i^{\bar a_i}\eta_{n_0}$ and denoting $\xi_{\bar a,\varepsilon}(x)=\xi_{\bar a}(\frac{x}{\varepsilon})$,  we have that $\xi_{\bar a}\in C^\infty_c({\R^N})$ and then for $\varepsilon\in(0,\, \frac12)$,
\begin{eqnarray*}
  \mathcal{L}^*_{\mu} \xi_{\bar a,\varepsilon}(x)  = \frac{1}{\varepsilon^2}  (-\Delta) \xi_{\bar a}(\frac{x}{\varepsilon})-\frac{1}{\varepsilon} \frac{x  }{|x|^2}  \cdot \nabla\xi_{\bar a}(\frac{x}{\varepsilon}),
\end{eqnarray*}
and on the one side,
\begin{eqnarray*}
\left| \int_{{\R^N}}u_{\mu}  \mathcal{L}^*_{\mu}( \xi_{\bar a,\varepsilon}) d\mu \right| &=& \left| \frac{1}{\varepsilon^2} \int_{B_{2\varepsilon}(0)}u_{\mu}\Gamma_{\mu} (-\Delta)  \xi_{\bar a}(\frac{x}{\varepsilon}) dx-  \frac{1}{\varepsilon}  \int_{B_{2\varepsilon}(0)}u_{\mu}\Gamma_{\mu} \frac{x  }{|x|^2}  \cdot \nabla\xi_{\bar a}(\frac{x}{\varepsilon})\right| \\
    &\le &\frac{1}{\varepsilon^2}\int_{B_{2\varepsilon}(0)}u_{\mu}\Gamma_{\mu} dx+\frac{1}{\varepsilon}     \int_{B_{2\varepsilon}(0)} \frac{ u_{\mu}\Gamma_{\mu} }{|x|} dx
    \\&\le&   \left\{\arraycolsep=1pt
\begin{array}{lll}
c_{12}  \quad
   &{\rm if}\quad N\ge 3,\\[1mm]
 \phantom{   }
- c_{12}\ln \varepsilon \quad  &{\rm   if}\quad N=2,
 \end{array}
 \right.
\end{eqnarray*}
where $c_{12}>0$ is independent of $\varepsilon$. Moreover, we have that
$$\left| \int_{{\R^N}}u_{\mu}   \xi_{\bar a,\varepsilon}   d\mu \right|\le \norm{\xi_{\bar a}}_{L^\infty} \int_{B_{2\varepsilon}(0)} f_{\mu}\Gamma_{\mu} dx\to0\quad{\rm as} \ \ \varepsilon\to 0^+. $$
On the other side,
$$\sum_{|a|=0}^p k_a D^{a} \xi_{\bar a,\varepsilon}(0) =\frac{k_{\bar a}}{\varepsilon^{|{\bar a}|}}|{\bar a!}|,  $$
where
$$|{\bar a}|=\sum\bar a_i\quad{\rm and}\quad \bar a!=\prod_{i=1}^N(\bar a_i)!\ge1.$$
So if $k_{\bar a}\not=0$, we have that
$$|\sum_{|a|=0}^p k_a D^{a} \xi_{\bar a,\varepsilon}(0)|\to +\infty\quad{\rm as} \ \ \varepsilon\to 0^+,$$
that is, the right hand of (\ref{4.2}) with $\xi=\xi_{\bar a,\varepsilon}$ blows up with the rate $\varepsilon^{-|\bar a|}$, the
which contradicts with the left hand of (\ref{4.2}) keeps bounded for $N\ge 3$ and blows up controlled by $-\ln \varepsilon$ as $\varepsilon\to0^+$,
so $k_a=0$ for $|a|\ge 1$.

Therefore, we have that
\begin{equation}\label{4.3}
 \int_{\R^N}(u_{\mu} \mathcal{L}^*_{\mu}\xi-  f\xi)\, d\mu  = k_0 \xi(0),\quad\forall\, \xi\in C^\infty_c(\R^N).
\end{equation}
For $\xi\in C^{1.1}_0(\R^N)$, take a sequence of functions in  $C^\infty_c(\R^N)$ converging to $\xi$, then the identity (\ref{4.3}) holds
for any $\xi\in C^{1.1}_c(\R^N)$.\hfill$\Box$

\medskip

\bigskip

\noindent{\bf Acknowledgements:} H. Chen  is supported by NNSF of China, No:11401270, 11661045, by the Jiangxi Provincial Natural Science Foundation, No: 20161ACB20007
and     by SRF for ROCS, SEM.


\begin{thebibliography}{99}



\bibitem{ACR} O. Adimurthi, N. Chaudhuri and M. Ramaswamy,
An improved Hardy-Sobolev inequality and its
application, {\it Proc. Amer. Math. Soc. 130}, 489-505 (2002).

\bibitem{A} P. Aviles, Local behaviour of the solutions of some elliptic equations,
{\it Comm. Math.
Phys. 108}, 177-192 (1987).

\bibitem {BOP}  L. Boccardo,  L. Orsina  and I. Peral,  A remark on existence and optimal summability of solutions of elliptic problems involving Hardy potential,
 {\it Discrete Contin. Dyn. Syst. 16},   513-523 (2006).



\bibitem {BDT} H. Brezis, L. Dupaigne and A. Tesei, On a semilinear elliptic equation with inverse-square potential,
{\it Selecta Mathematica 11.1}, 1-7 (2005).


\bibitem {BM} H. Brezis and  M. Marcus, Hardy's inequalities revisited, {\it Ann. Sc. Norm. Super. Pisa Cl. Sci. (5) 25},
 217-237 (1997).

\bibitem {BL} H. Brezis and P. Lions,
 A note on isolated singularities for linear elliptic equations, in Mathematical Analysis and Applications,  {\it Acad. Press},   263-266 (1981).

\bibitem {CW} F. Catrina  and Z. Wang,   On the Caffarelli-Kohn-Nirenberg inequalities: sharp constants, existence (and nonexistence), and symmetry of extremal functions,
 {\it Commun. Pure   Appl. Math., 54,} 229-258 (2001).

  \bibitem {CC} N. Chaudhuri and F. C\^irstea,
  On trichotomy of positive singular solutions associated with the Hardy-Sobolev operator,
  {\it Comptes Rendus Mathematique 347,}  153-158 (2009).

 \bibitem {CQZ} H. Chen, A. Quaas and F.  Zhou, On nonhomogeneous elliptic equations with the inverse square potential,
 {\it preprint.}

  \bibitem {DD} J. Davila and L. Dupaigne, Hardy-type inequalities, {\it J. Eur. Math. Soc. (JEMS) 6 (3)}, 335-365 (2004).


  \bibitem {D} L. Dupaigne,   A nonlinear elliptic PDE with the inverse square potential,
  {\it Journal d'Analyse Math\'ematique 86}, 359-398 (2002).

  \bibitem {F}  M. Fall, Nonexistence of distributional supersolutions of a semilinear elliptic equation with Hardy potential, {\it
   J. Funct. Anal. 264.3},  661-690 (2013).

 \bibitem {FF} V. Felli and A. Ferrero,  On semilinear elliptic equations with borderline Hardy potentials, {\it Journal D'Analyse Math\'ematique, 123,}
  303-340 (2014).

 \bibitem {GP} A. Garc\'ia  and G. Peral,   Hardy inequalities and some critical elliptic and parabolic problems,
 {\it  J. Differential Equations 144},  441-476 (1998).

\bibitem{GS} B. Gidas and J. Spruck, Global and local behaviour of positive solutions of nonlinear
elliptic equations, {\it Comm. Pure Appl. Math. 34}, 525-598 (1981).

\bibitem {G}  K. Gkikas,   Existence and nonexistence of energy solutions for linear elliptic equations involving Hardy-type potentials,
 {\it Indiana University Mathematics Journal 58}, 2317(2009).

 \bibitem {KV} K. Gkikas   and L. V\'eron, Boundary singularities of solutions of semilinear elliptic equations with critical Hardy potentials,
  {\it Nonlinear Analysis: Theory, Methods $\&$  Applications 121}, 469-540 (2015).

 \bibitem {IYM} K. Ishige,  K. Yoshitsu and O. El Maati,  The heat kernel of a Schr\" odinger operator with inverse square potential,
 {\it    arXiv:1602.04172} (2016).


\bibitem {MP} R. Mazzeo and F. Pacard,
A construction of singular solutions for a semilinear elliptic equation using asymptotic analysis,
{\it J. Differential Geom. 44},  331-370(1996).


 \bibitem{MN} M. Marcus and P. Nguyen,
 Moderate solutions of semilinear elliptic equations with Hardy potential,
 {\it  Ann. I. H. Poincar\'e -AN},  (2015).

  \bibitem {MT}
L. Moschini  and  A. Tesei, Parabolic Harnack inequality for the heat equation with inverse-square potential,
{\it Forum Mathematicum   19,} (2007).

\bibitem {P1} F. Pacard, Existence and convergence of positive weak solutions of $-\Delta u=u^{\frac{N}{N-2}}$ in bounded
domains of $\R^N$, {\it J. Calc. Var.  PDE. 1,} 243-265 (1993).



  \bibitem {P} A. Ponce, Elliptic PDEs, measures and capacities, {\it European Mathematical Society}, (2016).


  \bibitem {S} L. Schwartz, Theorie des distributions, {\it Hermann, Paris} (1966).

 \bibitem {V}  L. V\'{e}ron, Elliptic equations involving Measures,
 {\it Stationary Partial Differential equations,  M. Chipot, P.
 Quittner (Ed.) 1,}  593-712, (2004).
\end{thebibliography}
\end{document}